\documentclass[10pt]{amsart}
\usepackage{mathrsfs}
\usepackage{amsfonts}
\usepackage{amssymb}
\usepackage[papersize={6.7in,9.4in},textwidth=13.4cm,textheight=20cm,centering]{geometry}
\usepackage{enumerate}

\usepackage{multirow}

\usepackage{graphicx}
\usepackage{tikz}
\usepackage{float}

\usepackage[colorlinks=true,citecolor=blue,linkcolor=blue]{hyperref}
\hypersetup{
pdfstartpage=1,
pdfstartview=FitH}

\usepackage{amsthm,amssymb}
\usepackage{indentfirst}
\usepackage[leqno]{amsmath}

 \usepackage{caption} 
\numberwithin{equation}{section}

\allowdisplaybreaks

\newtheorem*{theorem*}{Theorem}
\newtheorem{theorem}{Theorem}[section]
\newtheorem{lemma}{Lemma}[section]

\newtheorem{corollary}{Corollary}[section]
\newtheorem{conjecture}{Conjecture}[section]
\newtheorem{theoremABC}{Theorem}[section]

\newtheorem*{claim*}{Claim}

\makeatletter
\newcounter{roem}
\renewcommand{\theroem}{\Roman{roem}}

\newcommand{\c@org@eq}{}
\let\c@org@eq\c@equation
\newcommand{\org@theeq}{}
\let\org@theeq\theequation

\newcommand{\setroem}{
\let\c@equation\c@roem
 \let\theequation\theroem}

\newcommand{\setarab}{
\let\c@equation\c@org@eq
\let\theequation\org@theeq}
\makeatother

\theoremstyle{definition}

\newcommand{\ue}{\mathrm{e}}
\newcommand{\ud}{\mathrm{d}}

\DeclareMathOperator{\Mod}{mod}

\renewcommand{\bmod}[1]{\,(\Mod{ #1})}

\newcommand{\bU}{\mathbf{U}}

\newcommand{\bR}{\mathbf{R}}
\newcommand{\bZ}{\mathbf{Z}}

\newcommand{\cD}{\mathcal{D}}

\newcommand{\cR}{\mathcal{R}}
\newcommand{\cS}{\mathcal{S}}

\newcommand{\fS}{\mathfrak{S}}

\newcommand{\rf}{\mathrm{f}}

\newcommand{\sB}{\mathscr{B}}
\newcommand{\sC}{\mathscr{C}}
\newcommand{\sD}{\mathscr{D}}

\newcommand{\sN}{\mathscr{N}}

\newcommand{\sU}{\mathscr{U}}
\newcommand{\sV}{\mathscr{V}}

\def\leq{\leqslant}

\begin{document}

\title[Counting fundamental solutions to Pell equation]
{Counting fundamental solutions to the Pell equation with prescribed size}

\author{Ping Xi}

\address{School of Mathematics and Statistics, Xi'an Jiaotong University, Xi'an 710049, P. R. China}

\email{ping.xi@xjtu.edu.cn}




\begin{abstract}
The cardinality of the set of $D\leqslant x$ for which
the fundamental solution of the Pell equation $t^2-Du^2=1$ is less than $D^{\frac{1}{2}+\alpha}$ with $\alpha\in[\frac{1}{2},1]$ is studied and certain lower bounds are obtained, improving previous results of Fouvry by introducing the $q$-analogue of van der Corput method to algebraic exponential sums with smooth moduli.
\end{abstract}

\maketitle


\section{Introduction and main results}\label{sec:Introduction}

Let $D$ be a non-square positive integer. The Pell equation is usually referred to
\begin{align}\label{eq:Pell}
t^2-Du^2=1,
\end{align}
to which the solution can be written in the usual form
\[\eta_D:=t+u\sqrt{D}.\]
The classical Dirichlet Units Theorem asserts that the set of solutions to \eqref{eq:Pell} is non-trivial and
has the form
\begin{align*}
\{\eta_D:\eta_D\text{ solution of }\eqref{eq:Pell}\}=\{\pm\varepsilon_D^n:n\in\bZ\},
\end{align*}
where $\varepsilon_D$ is called the fundamental solution of \eqref{eq:Pell} and is given by
\begin{align*}
\varepsilon_D:=\inf\{\eta_D:\eta_D>1\}.
\end{align*}
Writing $\varepsilon_D:=t_0+u_0\sqrt{D}$, we have $t_0,u_0\geqslant1$, from which we deduce that $t_0=\sqrt{1+u_0^2D}>\sqrt{D}$ and finally
\begin{align*}
\varepsilon_D>2\sqrt{D}.
\end{align*}
We are interested in counting the integers $D$ for which $\varepsilon_D$ or $\eta_D$ is less than a fixed power of $D$.

For $\alpha>0$ and $x\geqslant2$, define
\begin{align*}
S(x,\alpha)&:=|\{(\eta_D,D):2\leqslant D\leqslant x, ~D\neq\square,~\eta_D\leqslant D^{\frac{1}{2}+\alpha}\}|,\\
S^\rf(x,\alpha)&:=|\{(\varepsilon_D,D):2\leqslant D\leqslant x, ~D\neq\square,~\varepsilon_D\leqslant D^{\frac{1}{2}+\alpha}\}|.
\end{align*}
In his pioneer work, Hooley \cite{Ho84} proved the following theorem.
\begin{theoremABC} [Hooley]\label{thm:Hooley}
Let $\varepsilon_0$ satisfy $0<\varepsilon_0<\frac{1}{2}$. As $x\rightarrow+\infty$, one has
\begin{align*}
S(x,\alpha)=S^\rf(x,\alpha)=\Big(\frac{4\alpha^2}{\pi^2}+o(1)\Big)\sqrt{x}\log^2x
\end{align*}
uniformly for $\varepsilon_0\leqslant\alpha\leqslant\frac{1}{2}.$
\end{theoremABC}

In the same paper, Hooley \cite[P. 110]{Ho84} also made the following conjecture.
\begin{conjecture}[Hooley]\label{conjecture-Hooley}
For any given $\alpha>\frac{1}{2}$, we have
\begin{align*}
S^\rf(x,\alpha)=\frac{1}{\pi^2}\Big(4\alpha-1+C(\alpha)+o(1)\Big)\sqrt{x}\log^2x,
\end{align*}
where
\begin{align*}
C(\alpha)=\begin{cases}
0,\ \ &\alpha\in~]\frac{1}{2},1],\\
\frac{1}{18}(\alpha-1)^2,\ \ &\alpha\in~]1,\frac{5}{2}],\\
\frac{1}{24}(4\alpha-7),\ \ &\alpha\in~]\frac{5}{2},+\infty[.
\end{cases}
\end{align*}
\end{conjecture}

Fouvry \cite{Fo16} made a first significant step towards Hooley's conjecture in the case $\alpha\in~]\frac{1}{2},1].$ In fact, he proved the following theorem.
\begin{theoremABC} [Fouvry]\label{thm:Fouvry}
As $x\rightarrow+\infty$, one has
\begin{align}\label{eq:Fouvry-Sf}
S^\rf(x,\alpha)\geqslant \frac{1}{\pi^2}\Big(4\alpha-1-4\Big(\alpha-\frac{1}{2}\Big)^2-o(1)\Big)\sqrt{x}\log^2x
\end{align}
and
\begin{align}\label{eq:Fouvry-S}
S(x,\alpha)\geqslant\frac{1}{\pi^2}\Big(4\alpha-1-3\Big(\alpha-\frac{1}{2}\Big)^2-o(1)\Big)\sqrt{x}\log^2x
\end{align}
uniformly for $\alpha\in[\frac{1}{2},1].$
\end{theoremABC}

Moreover, Fouvry \cite{Fo16} considered a conjectural estimation for short exponential sums.
\begin{conjecture}\label{conjecture-expsums}
There exists an absolute $\vartheta\in[\frac{1}{2},1[$, 
such that for any integer $k\geqslant0,$ one has the inequality
\begin{align*}
\sum_{\substack{N<n\leqslant N_1\\ m\equiv a\bmod{4^k},(m,n)=1}}\ue\Big(\frac{h\overline{n}^2}{m^2}\Big)\ll_k(h,m^2)^{\frac{1}{2}}N^{\vartheta}
\end{align*}
uniformly for any integers $a, h, m$ satisfying $h\neq0,m\geqslant1$ and $2\nmid am$, for any real number $N$
satisfying $m\leqslant N\leqslant m^2$ and $N<N_1\leqslant2N.$
\end{conjecture}

Assuming Conjecture \ref{conjecture-expsums}, Fouvry \cite{Fo16} derived the following stronger lower bounds.
\begin{theoremABC} [Fouvry]\label{thm:Fouvry-conditional}
Assume that Conjecture $\ref{conjecture-expsums}$ is true for some $\vartheta\in[\frac{1}{2},1[$. As $x\rightarrow+\infty$, one has
\begin{align*}
S^\rf(x,\alpha)\geqslant \frac{1}{\pi^2}\Big(4\alpha-1-o(1)\Big)\sqrt{x}\log^2x
\end{align*}
and
\begin{align*}
S(x,\alpha)\geqslant\frac{1}{\pi^2}\Big(4\alpha-1+\Big(\alpha-\frac{1}{2}\Big)^2-o(1)\Big)\sqrt{x}\log^2x
\end{align*}
uniformly for $\alpha\in[\frac{1}{2},\frac{1}{1+\vartheta}].$
\end{theoremABC}
The first inequality coincides with Hooley's Conjecture \ref{conjecture-Hooley} when $\alpha$ is slightly larger than $\frac{1}{2}.$ On the other hand, Bourgain \cite{Bou15} considered Conjecture \ref{conjecture-expsums} itself. In particular, he succeeded in saving a power of $\log N$ in the trivial bound for the sum
involved. This allows him to improve the lower bound  \eqref{eq:Fouvry-Sf} by replacing the term $-4(\alpha-\frac{1}{2})^2$
with the term $O((\alpha-\frac{1}{2})^{2+c})$, where $c$ is some positive constant.

Bourgain's improvement is of interest only if $\alpha$ is quite close to $\frac{1}{2}$ and one should pay much more attention if every parameter is to be made effective. The aim of this paper is to give another improvement to Theorem \ref{thm:Fouvry} towards Conjecture \ref{conjecture-Hooley}. 
\begin{theorem}\label{thm:main}
Let $x\rightarrow+\infty$. For any fixed $\theta\in~]0,\frac{1}{2}[,$ we have
\begin{align}\label{eq:Sf(x,alpha)-lowerbound}
S^\rf(x,\alpha)
&\geqslant
\frac{1}{\pi^2}\Big\{4\alpha-1-4\Big(\alpha-\frac{1}{2}\Big)^2+\frac{1}{6}\rho\Big(\dfrac{1}{\theta}\Big)F_\theta(\alpha)-o(1)\Big\}\sqrt{x}\log^2x
\end{align}
and
\begin{align}\label{eq:S(x,alpha)-lowerbound}
S(x,\alpha)
&\geqslant
\frac{1}{\pi^2}\Big\{4\alpha-1-3\Big(\alpha-\frac{1}{2}\Big)^2+\frac{1}{6}\rho\Big(\dfrac{1}{\theta}\Big)F_\theta(\alpha)-o(1)\Big\}\sqrt{x}\log^2x
\end{align}
uniformly in $\alpha\in[\frac{1}{2},1]$, where $\rho$ is the Dickman function given by $\eqref{eq:Dickman}$ and
\begin{align}\label{eq:F}
F_\theta(\alpha)=
\begin{cases}
24\alpha-4(5+2\theta)\alpha^2-(7-2\theta),\ \ \ &\alpha\in[\frac{1}{2},\frac{6}{11+2\theta}],\\\noalign{\vskip 1mm}
864(11+2\theta)^{-2}-(7-2\theta),&\alpha\in~]\frac{6}{11+2\theta},1].
\end{cases}
\end{align}
\end{theorem}

Fix $\theta\in~]0,\frac{1}{2}[$. One may check that $F_\theta(\alpha)$ is always positive for
$\alpha\in~]\frac{1}{2},1]$ and monotonically increasing in
$\alpha\in[\frac{1}{2},\frac{6}{11+2\theta}]$. In particular, for $\alpha\in~]\frac{1}{2},\frac{6}{11}[,$ we take $\theta=\frac{3}{\alpha}-\frac{11}{2}$ such that
$\alpha=\frac{6}{11+2\theta}$
and
\begin{align*}
\frac{1}{6}\rho\Big(\dfrac{1}{\theta}\Big)F_\theta(\alpha)&=\rho\Big(\dfrac{2\alpha}{6-11\alpha}\Big)\frac{4(\alpha-\frac{1}{2})^3+6(\alpha-\frac{1}{2})^2}{\alpha}.
\end{align*}
Moreover, $\rho(\frac{2\alpha}{6-11\alpha})/\alpha>0.3008$ for $\alpha\in~]\frac{1}{2},\frac{35}{69}].$
Hence we may conclude the following consequence.
\begin{corollary}\label{corollary}
Let $x\rightarrow+\infty$. Uniformly for $\alpha\in~]\frac{1}{2},\frac{35}{69}],$ we have
\begin{align*}
S^\rf(x,\alpha)&\geqslant\frac{1}{\pi^2}\Big\{4\alpha-1-\frac{11}{5}\Big(\alpha-\frac{1}{2}\Big)^2+\frac{6}{5}\Big(\alpha-\frac{1}{2}\Big)^3\Big\}\sqrt{x}\log^2x
\end{align*}
and
\begin{align*}
S(x,\alpha)&\geqslant\frac{1}{\pi^2}\Big\{4\alpha-1-\frac{6}{5}\Big(\alpha-\frac{1}{2}\Big)^2+\frac{6}{5}\Big(\alpha-\frac{1}{2}\Big)^3\Big\}\sqrt{x}\log^2x.
\end{align*}
\end{corollary}

The framework of the proof is based on \cite{Fo16}. To make the paper clear, we will 
present the proof as complete as possible, but also with omitting some details that are not quite essential to understand the underlying ideas. 

The key point of proving Theorem \ref{thm:main} is a variant of Conjecture \ref{conjecture-expsums} that can be proved unconditionally. More precisely, if one allows the moduli $m^2$ to be smooth numbers (integers free of large prime factors), it is possible to prove the existence of $\vartheta$ in
Conjecture \ref{conjecture-expsums} as long as $N$ is not too small.
The details can be referred to Theorem \ref{thm:tripleexpsum-fSestimate} and Section \ref{sec:expsums}. We will adopt the $q$-analogue of van der Corput method, which can be at least dated back to Heath-Brown \cite{HB78} on the proof of Weyl-type subconvex bounds for Dirichlet $L$-functions to well-factorable moduli. Instead of the $AB$-process in \cite{HB78},
we apply the $BAB$-process by introducing a completion in the initial step.
It is expected that one can do better on the exponential sums in Conjecture \ref{conjecture-expsums} if better factorizations of the moduli are imposed; see 
\cite{WX16} for instance in the case of squarefree moduli. However, the improvements to Theorem \ref{thm:main} would be rather slight, since the density of smooth numbers decays rapidly when the size of their prime factors decreases.

As an extension to Theorem \ref{thm:main}, one may consider
$$S^\rf(x; \alpha, \beta):=|\{(\varepsilon_D,D): 2\leqslant D \leqslant x,~ D\neq\square,~ D^{\frac{1}{2} + \beta} \leqslant \varepsilon_D\leq D^{\frac{1}{2} + \alpha} \}|$$
for $\alpha>\beta\geqslant0.$
Conjecture \ref{conjecture-Hooley} would yield asymptotics for $S(x; \alpha, \beta)$ while $\alpha,\beta$ are of different prescribed sizes.
A weaker statement would assert that, for any $\alpha>\beta\geqslant0,$ there exists a positive constant $c=c(\alpha,\beta)$, such that 
$$S^\rf(x; \alpha, \beta)\geqslant c\sqrt{x}\log^2x$$ 
for all large $x>x_0(\alpha,\beta)$.
This weaker statement was made unconditionally by Fouvry and Jouve \cite{FJ12}
whenever $\beta<\frac{3}{2}$. It is expected that the arguments in this paper can enlarge the admissible range of $\beta$.

\subsection*{Notation and convention} 
As usual, $\ue(x)=\ue^{2\pi ix}$, $\varphi$ denotes the Euler function and $\omega(q)$ counts the number of distinct prime factors of $q$. The variable $p$ is reserved for prime numbers. Denote by $q^\flat$ and $q^\sharp$ the squarefree and squarefull parts of $q$, respectively; namely,
\[q^\flat=\prod_{p\| q}p,\ \ \ q^\sharp=\prod_{p^\nu\|q,\nu\geqslant2}p^\nu.\]
For a real number $x,$ denote by $[x]$ its integral part and $\|x\|=\min_{n\in\bZ}|x-n|$. 
From time to time, we use $(m,n)$ to denote the greatest common divisor of $m,n$, and also to denote a tuple given by two coordinates; these will not cause confusions
as one will see later. The symbol $*$ in summation reminds us to sum over primitive elements such that poles of the summand are avoided.
For a function $f\in L^1(\bR)$, its Fourier transform is defined as
\begin{align*}
\widehat{f}(y) := \int_\bR f(x) \mathrm{e}(-yx)\ud x.
\end{align*}

We use $\varepsilon$ to denote a very small positive number, which might be different at each occurrence; we also write $X^\varepsilon \log X\ll X^\varepsilon.$ 
The convention $n\sim N$ means $N<n\leqslant2N.$
\subsection*{Acknowledgements} 
I am very grateful to the referee for the valuable comments and suggestions.
The work is supported in part by NSFC (No.11601413) and NSBRP (No. 2017JQ1016) of Shaanxi Province.

\smallskip

\section{Fundamental transformations: after Hooley and Fouvry}
We first make some fundamental transformations following the arguments of Hooley and Fouvry. For some conclusions, we omit the proof and the detailed arguments can be found in \cite{Ho84} and \cite{Fo16}.

\subsection{An initial transformation}
First, we write
\begin{align}\label{eq:reduction-initial}
S(x,\alpha)&=\mathop{\sum_{\square\neq D\leqslant x}\sum_{t\geqslant1}\sum_{u\geqslant1}}_{\substack{t^2-Du^2=1\\1<t+u\sqrt{D}\leqslant D^{\frac{1}{2}+\alpha}}}1=\sum_{1\leqslant u\leqslant X_\alpha}\sum_{\substack{Y_2(u,\alpha)\leqslant t\leqslant Y_3(u)\\ Du^2=t^2-1}}1,
\end{align}
where
\begin{align}
X_\alpha&:=\frac12(x^\alpha-x^{-1-\alpha});\label{eq:X_alpha}\\
Y_3(u)&:=\sqrt{xu^2+1};\label{eq:Y_3(u)}\\
Y_2(u,\alpha)&:=\sqrt{Y_1(u,\alpha)u^2+1}.\label{eq:Y_2(u,alpha)}
\end{align}
Here, $Y_1(u,\alpha)$ is a function in $u$, implicitly defined by the equation
\begin{align*}
u=\frac12(Y_1(u,\alpha)^\alpha-Y_1(u,\alpha)^{-1-\alpha}).
\end{align*}

We have the following asymptotic characterization for $Y_2(u,\alpha).$ The proof can be found in \cite[Lemma 2.1]{Fo16} and the subsequent remark.
\begin{lemma}\label{lm:Y_2(u,alpha)}
Let $\alpha>0.$ The function $u\mapsto Y_2(u,\alpha)$ is of $\sC^\infty$-class and satisfies the inequalities
\begin{align}\label{eq:Y_2(u)-1}
2^{\frac{1}{2\alpha}}u^{1+\frac{1}{2\alpha}}<Y_2(u,\alpha)<(2^{\frac{1}{2\alpha}}+o(1))u^{1+\frac{1}{2\alpha}}
\end{align}
and
\begin{align}\label{eq:Y_2(u)-2}
\frac{\ud }{\ud u}Y_2(u,\alpha)=O(u^{\frac{1}{2\alpha}})
\end{align}
as $u\rightarrow+\infty.$
\end{lemma}

\subsection{A first dissection of $S(x,\alpha)$}
We truncate the $u$-sum in \eqref{eq:reduction-initial} at $X_{\frac{1}{2}}$, and the 
contributions from $u\leqslant X_{\frac{1}{2}}$ and $u>X_{\frac{1}{2}}$
are respectively denoted by $A(x,\alpha)$ and $B(x,\alpha).$ Therefore,
\begin{align}\label{eq:S=A+B}
S(x,\alpha)=A(x,\alpha)+B(x,\alpha).
\end{align}
Accordingly, we may define $A^\rf(x,\alpha)$ and $B^\rf(x,\alpha)$
by introducing the extra restriction $t+u\sqrt{D}=\varepsilon_D$
to the equation $t^2-Du^2=1.$ 

As stated by Fouvry \cite[Formula (4.5)]{Fo16}, we have
\begin{lemma} Let $\alpha\in[\frac{1}{2},1].$ As $x\rightarrow+\infty,$ we have
\begin{align}\label{eq:A(x,alpha)-asymptotic}
A(x,\alpha)=\Big(\frac{1}{\pi^2}+o(1)\Big)\sqrt{x}\log^2x.
\end{align}
\end{lemma}

Our task thus reduces to proving a lower bound for $B(x,\alpha).$ Put
$\cR(u):=\{\Omega\bmod{u^2}:\Omega^2\equiv1\bmod{u^2}\}.$ We then have
\begin{align}\label{eq:B(x,alpha)}
B(x,\alpha)=\sum_{X_{\frac{1}{2}}< u\leqslant X_\alpha}\sum_{\Omega\in\cR(u)}\sum_{\substack{Y_2(u,\alpha)\leqslant t\leqslant Y_3(u)\\ t\equiv\Omega\bmod{u^2}}}1.
\end{align}

Put $\gamma(u)=|\cR(u)|$. Thus $u\mapsto \gamma(u)$
is multiplicative and satisfies
\begin{align}\label{eq:gamma-primepower}
\begin{cases}
\gamma(2)=2, \gamma(2^k)=4\text{ for }k\geqslant2,\\
\gamma(p^\ell)=2\text{ for }p\geqslant3,\ell\geqslant1.
\end{cases}
\end{align}

\subsection{Analysis of $\cR(u)$}
For $u\geqslant1$, write $u=2^ku_0$, where $u_0$ is an odd integer. 
The choice of $(k,u_0)$ is unique for each $u\geqslant1$.
The Chinese remainder theorem implies
\begin{align}\label{eq:R(u)-equivalence}
\Omega^2\equiv1\bmod{u^2}\Longleftrightarrow
\begin{cases}
\Omega^2\equiv1\bmod{u_0^2},\\
\Omega^2\equiv1\bmod{4^k}.\end{cases}
\end{align}
In this way, one can establish a bijection between $\cR(u)$ and $\cR(2^k)\times\cR(u_0).$
Starting from \eqref{eq:B(x,alpha)}, we decompose $B(x,\alpha)$ by
\begin{align}\label{eq:B(x,alpha)-decomposition}
B(x,\alpha)=\sum_{k\geqslant0}\sum_{\xi\in\cR(2^k)}B(x,\alpha;\xi,k),
\end{align}
where
\begin{align*}
B(x,\alpha;\xi,k)=\sum_{\substack{2^{-k}X_{\frac{1}{2}}< u\leqslant 2^{-k}X_\alpha\\2\nmid u}}\sum_{\Omega\in\cR(u)}\sum_{\substack{Y_2(2^ku,\alpha)\leqslant t\leqslant Y_3(2^ku)\\ t\equiv\Omega\bmod{u^2}\\ t\equiv\xi\bmod{4^k}}}1.
\end{align*}

The task will be evaluating $B(x,\alpha;\xi,k)$
for all $k\geqslant0$ and $\xi\in\cR(2^k).$ This would require the following description of $\cR(u)$ that allows us to create one more variable. This is Lemma 4.1 in \cite{Fo16}.

\begin{lemma}\label{lm:R(u)-description}
Let $u$ be a positive odd integer. Then there is a bijection $\Phi$ between the set
of coprime decompositions of $u$
\[\cD(u):=\{(u_1,u_2):u_1u_2=u,~(u_1,u_2)=1,~u_1,u_2\geqslant1\}\]
and the set of roots of congruence
\[\cR(u):=\{\Omega\bmod{u^2}:\Omega^2\equiv1\bmod{u^2}\}.\]
Such a bijection can be defined by $\Phi(u_1,u_2)=\Omega,$ where $\Omega$ is uniquely determined by the congruences
$\Omega\equiv1\bmod{u_1^2}$ and $\Omega\equiv-1\bmod{u_2^2}$ .
In an equivalent manner, we have the congruence
\begin{align}\label{eq:Phi-congruence}
\Phi(u_1,u_2)&\equiv-\overline{u}_1^2u_1^2+\overline{u}_2^2u_2^2\bmod{u^2}.
\end{align}
Here $\overline{u}_1u_1\equiv1\bmod{u_2}$ and $\overline{u}_2u_2\equiv1\bmod{u_1}$.
\end{lemma}

With the help of Lemma \ref{lm:R(u)-description}, we may rewrite
$B(x,\alpha;k,\xi)$ as
\begin{align*}
B(x,\alpha;\xi,k)&=\mathop{\sum\sum}_{\substack{2^{-k}X_{\frac{1}{2}}< u_1u_2\leqslant 2^{-k}X_\alpha\\(2,u_1u_2)=(u_1,u_2)=1}}\sum_{\substack{Y_2(2^ku_1u_2,\alpha)\leqslant t\leqslant Y_3(2^ku_1u_2)\\ t\equiv\Phi(u_1,u_2)\bmod{u_1^2u_2^2}\\ t\equiv\xi\bmod{4^k}}}1\\
&=:B^>(x,\alpha;\xi,k)+B^<(x,\alpha;\xi,k),
\end{align*}
where $B^>(x,\alpha;\xi,k)$ and $B^<(x,\alpha;\xi,k)$ correspond to the restrictions
$u_1>u_2$ and $u_1<u_2$, respectively.
Since the treatments of $B^>(x,\alpha;\xi,k)$ and $B^<(x,\alpha;\xi,k)$
are similar, it suffices to study $B^<(x,\alpha;\xi,k)$ as presented in the next section.

We close this section by the following trivial equality:
\begin{align}\label{eq:S(x,alpha)-Sf(x,alpha)}
S(x,\alpha)=S^\rf(x,\alpha)+S\Big(x,\frac{\alpha}{2}-\frac{1}{4}\Big),\ \ \frac{1}{2}<\alpha\leqslant1,
\end{align}
which is a consequence of the equivalence
\[\varepsilon_D^2\leqslant D^{\frac{1}{2}+\alpha}\Longleftrightarrow \varepsilon_D\leqslant D^{\frac{1}{2}+(\frac{\alpha}{2}-\frac{1}{4})}.\]
This allows us to transfer between $S^\rf(x,\alpha)$ and $S(x,\alpha)$.

\smallskip

\section{Lower bound for $B(x,\alpha)$}
In order to conclude the lower bound for $B(x,\alpha)$,
we now start the study of $B^<(x,\alpha;\xi,k)$.
Recall that
\begin{align*}
B^<(x,\alpha;\xi,k)&=\mathop{\sum\sum}_{\substack{2^{-k}X_{\frac{1}{2}}< u_1u_2\leqslant 2^{-k}X_\alpha\\u_1<u_2\\(2,u_1u_2)=(u_1,u_2)=1}}\sum_{\substack{Y_2(2^ku_1u_2,\alpha)\leqslant t\leqslant Y_3(2^ku_1u_2)\\ t\equiv\Phi(u_1,u_2)\bmod{u_1^2u_2^2}\\ t\equiv\xi\bmod{4^k}}}1.
\end{align*}
We would like to drop the multiplicative constraints $2^{-k}X_{\frac{1}{2}}< u_1u_2\leqslant 2^{-k}X_\alpha$ and sum over $u_1,u_2$ separately.
To do so, we may introduce the following inequality
\begin{align}\label{eq:B^<-lowerbound}
B^<(x,\alpha;\xi,k)&\geqslant\sum_{\bU}\sum_{\xi_1}\sum_{\xi_2}B(x,\alpha;\bU,\xi,\xi_1,\xi_2,k),
\end{align}
where 
\begin{itemize}
\item the summation is over all $\bU=(U_1,U_2)$ satisfying
\begin{align*}
U_1<U_2,\ \ \ X_{\frac{1}{2}}< 2^kU_1U_2\leqslant \frac{X_\alpha}{8},
\end{align*}
and $U_1,U_2$ being powers of $2$,
\item the summation is over all $\xi_1,\xi_2\bmod{4^k}$ satisfying $(\xi_1\xi_2,4^k)=1,$
\item we have defined
\begin{align*}
B(x,\alpha;\bU,\xi,\xi_1,\xi_2,k)&:=
\mathop{\sum\sum}_{\substack{u_1\sim U_1,u_2\sim U_2\\u_1\equiv\xi_1,u_2\equiv\xi_2\bmod{4^k}\\(2,u_1u_2)=(u_1,u_2)=1}}\sum_{\substack{Y_2(2^ku_1u_2,\alpha)\leqslant t\leqslant Y_3(2^ku_1u_2)\\ t\equiv\Phi(u_1,u_2)\bmod{u_1^2u_2^2}\\ t\equiv\xi\bmod{4^k}}}1.
\end{align*}
\end{itemize}
Of course the condition $(2,u_1u_2)=1$ can be dropped when $k\geqslant1.$
The parameter $\alpha$ is supposed
to be fixed and the congruence conditions modulo $4^k$ are harmless. So to shorten the
notations, we write $B(x,\bU):=B(x,\alpha;\bU,\xi,\xi_1,\xi_2,k).$
Finally we shall not precise the dependence on $k$ of some $O$-symbols, since we shall work with
a finite number of values of $k$. The case $k=0$ is typical
and really reflects the difficulties of the method.

\subsection{Reduction to exponential sums: after Fouvry}
The congruence condition $t\equiv\Phi(u_1,u_2)\bmod{u_1^2u_2^2}$
implies that $t\equiv-1\bmod{u_2^2}$, i.e., $t=-1+\ell u_2^2$ for some $\ell\in\bZ.$
Since $Y_2(2^ku_1u_2,\alpha)\leqslant t\leqslant Y_3(2^ku_1u_2)$, then there is no such $t$ if $u_2$ is too large, for instance when
\[u_2^2>\sqrt{4^kxu_1^2u_2^2+1}+1.\]
Hence we can suppose
\begin{align}\label{eq:U1U2-sizes}
U_2\leqslant 2^{k+2}\sqrt{x}U_1,
\end{align}
otherwise $B(x,\bU)=0.$

Since $u_1u_2$ is odd, we deduce from \eqref{eq:Phi-congruence} the equivalence
\begin{align*}
t\equiv\Phi(u_1,u_2)\bmod{u_1^2u_2^2},\ t\equiv\xi\bmod{4^k}\Longleftrightarrow t\equiv t_0\bmod{4^ku_1^2u_2^2}
\end{align*}
with
\begin{align*}
t_0\equiv\xi u_1^2u_2^2\cdot\overline{(u_1^2u_2^2)}
+(2u_2^2\overline{u_2^2}-1)4^k\overline{4^k}\bmod{4^ku_1^2u_2^2},
\end{align*}
where
$u_1^2u_2^2\cdot\overline{(u_1^2u_2^2)}\equiv1\bmod{4^k},
u_2^2\overline{u_2^2}\equiv1\bmod{u_1^2}$ and $4^k\overline{4^k}\equiv1\bmod{u_1^2u_2^2}.$
It follows that
\begin{align*}
\frac{t_0}{4^ku_1^2u_2^2}\equiv\kappa-\frac{1}{4^ku_1^2u_2^2}+2\frac{\overline{4^ku_2^2}}{u_1^2}\bmod1
\end{align*}
with $\kappa:=(\xi+1)\overline{\xi_1^2\xi_2^2}/4^k.$
The three terms on the RHS have completely different structures: the first one
is constant, the second one changes very slowly when $u_1$ and $u_2$ vary, the third one oscillates
a lot when $u_2$ varies with $u_1$ fixed. 

For each fixed $k$, we rewrite the sum $B(x,\bU)$ as
\begin{align*}
B(x,\bU)=B(x,\alpha;\bU,\xi,\xi_1,\xi_2,k)=
\mathop{\sum\sum}_{\substack{u_1\sim U_1,u_2\sim U_2\\u_1\equiv\xi_1,u_2\equiv\xi_2\bmod{4^k}\\(2,u_1u_2)=(u_1,u_2)=1}}
\sum_{\substack{Y_2\leqslant t\leqslant Y_3\\ t\equiv t_0\bmod{4^ku_1^2u_2^2}}}1,
\end{align*}
where $Y_2:=Y_2(2^ku_1u_2,\alpha)$ and $Y_3:=Y_3(2^ku_1u_2)$.
As in \cite{Fo16}, we smooth the $t$-sum via the following lemma.
\begin{lemma}\label{lm:smoothfunction}
For every $\delta>0$ there exists a smooth function $g:\bR\rightarrow\bR$ which has the two
properties
\[0\leqslant g\leqslant \mathbf{1}_{[-\frac{1}{2},\frac{1}{2}]}\]
and
\[\int_\bR g(y)\ud y=1-\delta.\]
\end{lemma}

Let $g$ be a smooth function given as in Lemma \ref{lm:smoothfunction}. Hence
\begin{align}\label{eq:B(x,bU)-lowerbound}
B(x,\bU)&\geqslant
\mathop{\sum\sum}_{\substack{u_1\sim U_1,u_2\sim U_2\\u_1\equiv\xi_1,u_2\equiv\xi_2\bmod{4^k}\\(2,u_1u_2)=(u_1,u_2)=1}}
\sum_{\substack{t\in\bZ\\ t\equiv t_0\bmod{4^ku_1^2u_2^2}}}g\Big(\frac{t-\frac{Y_2+Y_3}{2}}{Y_3-Y_2}\Big).
\end{align}
By Poisson summation, the $t$-sum becomes
\begin{align*}
\frac{Y_3-Y_2}{4^ku_1^2u_2^2}\sum_{h\in\bZ}\widehat{g}\Big(\frac{h(Y_3-Y_2)}{4^ku_1^2u_2^2}\Big)\ue\Big(h\kappa+2h\frac{\overline{4^ku_2^2}}{u_1^2}-\frac{h(2+Y_2+Y_3)}{2\cdot 4^ku_1^2u_2^2}\Big).
\end{align*}
From integration by parts, we have $\widehat{g}(y)\ll(1+|y|)^{-A}$ for any $A\geqslant0$. Note that
\begin{align*}
\frac{Y_3-Y_2}{4^ku_1^2u_2^2}\asymp\frac{\sqrt{x}}{U_1U_2}.
\end{align*}
The above sum over $h$ can be truncated to $0\leqslant|h|\leqslant H$ with $H=U_1U_2x^{-\frac{1}{2}+\varepsilon}$, and the remaining contribution is at most $O(x^{-2017}).$
Therefore,
\begin{align*}
B(x,\bU)&\geqslant
\mathop{\sum\sum}_{\substack{u_1\sim U_1,u_2\sim U_2\\u_1\equiv\xi_1,u_2\equiv\xi_2\bmod{4^k}\\(2,u_1u_2)=(u_1,u_2)=1}}
\frac{Y_3-Y_2}{4^ku_1^2u_2^2}\sum_{0\leqslant|h|\leqslant H}\widehat{g}\Big(\frac{h(Y_3-Y_2)}{4^ku_1^2u_2^2}\Big)\\
&\ \ \ \ \times\ue\Big(h\kappa+2h\frac{\overline{4^ku_2^2}}{u_1^2}-\frac{h(2+Y_2+Y_3)}{2\cdot 4^ku_1^2u_2^2}\Big)+O(x^{-1949})\\
&=B_1(x,\bU)+B_2(x,\bU)+O(x^{-1949}),
\end{align*}
where $B_1(x,\bU)$ and $B_2(x,\bU)$ are used to denote contributions from $h=0$
and $h\neq0$, respectively.

First,
\begin{align*}
B_1(x,\bU)&=(1-\delta)\mathop{\sum\sum}_{\substack{u_1\sim U_1,u_2\sim U_2\\u_1\equiv\xi_1,u_2\equiv\xi_2\bmod{4^k}\\(2,u_1u_2)=(u_1,u_2)=1}}
\frac{Y_3-Y_2}{4^ku_1^2u_2^2},
\end{align*}
which is $\asymp\sqrt{x}$. It is also desirable to show that
\begin{align}\label{eq:B2(x,U)-upperbound}
B_2(x,\bU)\ll x^{\frac{1}{2}-\varepsilon_0}
\end{align}
for some $\varepsilon_0>0$. By standard tools from analysis (see \cite{Fo16} for details), it suffices to prove that
\begin{align}\label{eq:tripleexpsum-expected}
\sum_{\substack{h\leqslant H\\ h\equiv h_0\bmod{4^k}}}\mathop{\sum_{u_1\in]U_1,U_1^*]}\sum_{u_2\in]U_2,U_2^*]}}
_{\substack{u_1\equiv\xi_1,u_2\equiv\xi_2\bmod{4^k}\\(2,u_1u_2)=(u_1,u_2)=1}}
\ue\Big(2h\frac{\overline{4^ku_2^2}}{u_1^2}\Big)\ll U_1U_2x^{-\varepsilon_0},
\end{align}
where $U_j<U_j^*\leqslant2U_j$, $j=1,2$. After transforming the $u_2$-sum to a complete sum $T(\cdot,u_1^2)$, where 
\begin{align*}
T(h,q)=\sideset{}{^*}\sum_{x\bmod q}\ue\Big(\frac{h\overline{x^2}}{q}\Big)=\sideset{}{^*}\sum_{x\bmod q}\ue\Big(\frac{hx^2}{q}\Big),
\end{align*}
Fouvry \cite{Fo16} evaluated $T(\cdot,u_1^2)$ in terms of classical Gauss sums and Jacobi symbols. He then arrived at a bilinear form involving Jacobi symbols, for which a celebrated estimate due to Heath-Brown \cite{HB95} was applied.
Amongst some other delicate arguments, Fouvry was able to prove
\eqref{eq:B2(x,U)-upperbound} under the conditions
\begin{align}\label{eq:U1U2-sizes-Fouvry}
U_1\leqslant x^{\frac{1}{4}-5\varepsilon_0},\ \ \ U_2\leqslant U_1 x^{\frac{1}{2}-5\varepsilon_0},
\end{align}
in which case he obtained the lower bound
\begin{align}\label{eq:B(x,U)-lowerbound-Fouvry}
B(x,\bU)\geqslant (1-\delta)\mathop{\sum\sum}_{\substack{u_1\sim U_1,u_2\sim U_2\\u_1\equiv\xi_1,u_2\equiv\xi_2\bmod{4^k}\\(2,u_1u_2)=(u_1,u_2)=1}}
\frac{Y_3-Y_2}{4^ku_1^2u_2^2}+O(x^{\frac{1}{2}-\varepsilon_0}).
\end{align}

To obtain a better lower bound for $B(x,\alpha)$ and thus for $S(x,\alpha)$, it is natural to expect that \eqref{eq:B2(x,U)-upperbound} and \eqref{eq:B(x,U)-lowerbound-Fouvry} can hold in larger ranges of $U_1,U_2$.
However, it seems rather difficult when $U_1$ is quite close to $U_2$ since the $u_2$-sum is too short in the sense of the P\'olya-Vinogradov barrier. In fact, Bourgain \cite{Bou15} managed to control the LHS in 
\eqref{eq:tripleexpsum-expected}, but with a saving of a small power of $\log x$ rather than that of $x$. This allows him to improve upon Fouvry when $\alpha$ is rather close to $\frac{1}{2}$ in Theorem \ref{thm:Fouvry}.
 
In our subsequent argument, we will specialize $u_1$ with special structures in the original sum \eqref{eq:B(x,bU)-lowerbound} before Poisson summation.
More precisely, we will consider those $u_1$ consisting of only small prime factors, so that 
$u_1$ has good factorizations, which enable us to control the exponential sums in
\eqref{eq:tripleexpsum-expected} even though $U_1$ is quite close to $U_2$.

\subsection{Lower bound of $B(x,\bU)$: smooth approach}
A positive integer $n$ is said to be $y$-smooth (or friable) if all prime factors of $n$
do not exceed $y$. Let $\theta\in~]0,\frac{1}{2}[$ be a fixed number.
If $n$ is $n^{\theta}$-smooth, the inclusion-exclusion principle yields the existence
of the divisor $d\mid n$ such that
$n^{\theta_0}\leqslant d\leqslant n^{\theta_0+\theta}$
for any $\theta_0\in[0,1-\theta].$

We now restrict these $u_1$ in the RHS of 
\eqref{eq:B(x,bU)-lowerbound} to $U_1^{\theta}$-smooth numbers
and put
\begin{align}\label{eq:B*(x,bU)-lowerbound}
B^*(x,\bU)&=\mathop{\sum\sum}_{\substack{u_1\sim U_1,u_2\geqslant1\\u_1\equiv\xi_1,u_2\equiv\xi_2\bmod{4^k}\\(2,u_1u_2)=(u_1,u_2)=1\\ u_1\text{ is $U_1^{\theta}$-smooth}}}
g_1\Big(\frac{u_2}{U_2}\Big)\sum_{\substack{t\in\bZ\\ t\equiv t_0\bmod{4^ku_1^2u_2^2}}}g\Big(\frac{t-\frac{Y_2+Y_3}{2}}{Y_3-Y_2}\Big),
\end{align}
where $g_1(y)=g(y-\frac{3}{2})$ with $g$ given as in Lemma \ref{lm:smoothfunction}.
Following the similar arguments of Fouvry, we may derive that
\begin{align*}
B^*(x,\bU)&\geqslant B_1^*(x,\bU)+B_2^*(x,\bU)+O(x^{-1949}),
\end{align*}
where 
\begin{align*}
B_1^*(x,\bU)
&=(1-\delta)\mathop{\sum\sum}_{\substack{u_1\sim U_1,u_2\geqslant1\\u_1\equiv\xi_1,u_2\equiv\xi_2\bmod{4^k}\\(2,u_1u_2)=(u_1,u_2)=1\\ u_1\text{ is $U_1^{\theta}$-smooth}}}g_1\Big(\frac{u_2}{U_2}\Big)
\frac{Y_3-Y_2}{4^ku_1^2u_2^2}\\
&\geqslant(1-\delta)^2\mathop{\sum\sum}_{\substack{u_1\sim U_1,u_2\sim U_2\\u_1\equiv\xi_1,u_2\equiv\xi_2\bmod{4^k}\\(2,u_1u_2)=(u_1,u_2)=1\\ u_1\text{ is $U_1^{\theta}$-smooth}}}
\frac{Y_3-Y_2}{4^ku_1^2u_2^2},
\end{align*}
and we expect to show that
\begin{align}\label{eq:B2*(x,U)-upperbound}
B_2^*(x,\bU)\ll x^{\frac{1}{2}-\varepsilon_0}
\end{align}
for some $\varepsilon_0>0$, for which it suffices to prove, for all $U_1^*\in~]U_1,2U_1]$, that
\begin{align}\label{eq:tripleexpsum-expected-smooth}
\sum_{\substack{h\leqslant H\\ h\equiv h_0\bmod{4^k}}}\mathop{\sum_{u_1\in]U_1,U_1^*]}\sum_{u_2\geqslant1}}
_{\substack{u_1\equiv\xi_1,u_2\equiv\xi_2\bmod{4^k}\\(2,u_1u_2)=(u_1,u_2)=1\\ u_1\text{ is $U_1^{\theta}$-smooth}}}
g_1\Big(\frac{u_2}{U_2}\Big)\ue\Big(2h\frac{\overline{4^ku_2^2}}{u_1^2}\Big)\ll U_1U_2x^{-\varepsilon_0}.
\end{align}

We will prove
\begin{theorem}\label{thm:tripleexpsum-fSestimate}
There exists some $\varepsilon_0\in~]0,10^{-2017}[$ such that $\eqref{eq:tripleexpsum-expected-smooth}$ holds, provided that
\begin{align}\label{eq:U1U2-sizes-smooth}
U_1^{\frac{8+2\theta}{3}}U_2\leqslant x^{1-2\varepsilon_0},\ \ \ U_2\leqslant U_1 x^{\frac{1}{2}-5\varepsilon_0}.
\end{align}
\end{theorem}
Put $U_1=x^{\gamma_1}$ and $U_2=x^{\gamma_2}$.
In addition to the restrictions $\gamma_2<\gamma_1+\frac{1}{2},\frac{1}{2}<\gamma_1+\gamma_2<\alpha$, \eqref{eq:U1U2-sizes-Fouvry} requires $\gamma_1<\frac{1}{4}$, and 
we require
$\gamma_2<1-\frac{26}{9}\gamma_1$
in the particular case $\theta=\frac{1}{36}$. In the following figure, the blue area shows what we can gain more than the previous approach (We are gaining relatively more as
$\alpha$ becomes closer to $\frac{1}{2}$).
\begin{figure}[H]
\includegraphics[width=8.5cm]{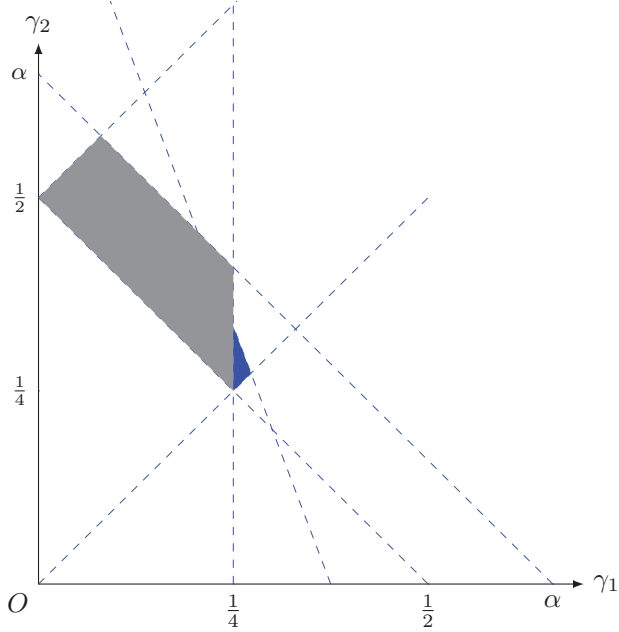}
\caption{Admissible values of $(\gamma_1,\gamma_2)$}
\end{figure}

The proof of Theorem \ref{thm:tripleexpsum-fSestimate} will be given in the next section.
To see the advantage of our approach, one may consider the particular case $U_1=U_2$, 
and our first restriction will reduce to $U_1\leqslant x^{\frac{3}{11+2\theta}-\varepsilon_0}$; however, the stronger restriction $U_1\leqslant x^{\frac{1}{4}-5\varepsilon_0}$ in \eqref{eq:U1U2-sizes-Fouvry} is required.

Therefore, 
we may obtain the lower bound 
\begin{align}\label{eq:B(x,U)-lowerbound-smooth}
B(x,\bU)&\geqslant(1-\delta)^2\mathop{\sum\sum}_{\substack{u_1\sim U_1,u_2\sim U_2\\u_1\equiv\xi_1,u_2\equiv\xi_2\bmod{4^k}\\(2,u_1u_2)=(u_1,u_2)=1\\ u_1\text{ is $U_1^{\theta}$-smooth}}}
\frac{Y_3-Y_2}{4^ku_1^2u_2^2}+O(x^{\frac{1}{2}-\varepsilon_0})
\end{align}
subject to the restrictions in \eqref{eq:U1U2-sizes-smooth}.

\subsection{A weakened form of Theorem \ref{thm:main}}
Up to now, we have two lower bounds for $B(x,\bU)$, i.e.,
\eqref{eq:B(x,U)-lowerbound-Fouvry} and \eqref{eq:B(x,U)-lowerbound-smooth},
subject to the restrictions in \eqref{eq:U1U2-sizes-Fouvry} and \eqref{eq:U1U2-sizes-smooth}, respectively.
In what follows, we will take into account all such admissible tuples $(k,\xi,\xi_1,\xi_2,U_1,U_2)$,
for which we appeal to \eqref{eq:B(x,U)-lowerbound-Fouvry} if \eqref{eq:U1U2-sizes-Fouvry} is satisfied, and appeal to \eqref{eq:B(x,U)-lowerbound-smooth} if 
\eqref{eq:U1U2-sizes-Fouvry} is not satisfied but \eqref{eq:U1U2-sizes-smooth} is valid.
To this end, we define two sets of tuples
\begin{align*}
\sU_1:=\{(U_1,U_2):U_1<U_2,\ \frac{\sqrt{x}}{2}< 2^kU_1U_2\leqslant \frac{x^\alpha}{16},\  U_1\leqslant x^{\frac{1}{4}-\eta},\ U_2\leqslant U_1 x^{\frac{1}{2}-\eta}\},
\end{align*}
and
\begin{align*}
\sU_2:=\{(U_1,U_2):U_1<U_2,\ \frac{\sqrt{x}}{2}< 2^kU_1U_2\leqslant \frac{x^\alpha}{16},\  U_1^{\frac{8+2\theta}{3}}U_2\leqslant x^{1-\eta},\ U_2\leqslant U_1 x^{\frac{1}{2}-\eta}\}\setminus \sU_1,
\end{align*}
where $\eta$ is a sufficiently small positive number.

First, we may derive a lower bound for $B^<(x,\alpha;\xi,k)$ by 
inserting the inequality \eqref{eq:B(x,U)-lowerbound-Fouvry} or \eqref{eq:B(x,U)-lowerbound-smooth} to \eqref{eq:B^<-lowerbound}. A similar lower bound also holds for $B^>(x,\alpha;\xi,k)$ by symmetry.
Therefore, we have
\begin{align*}
B(x,\alpha)&\geqslant2(1-\delta)\sum_{k=0}^{k_0}\sum_{\xi\in\cR(2^k)}\sum_{\bU\in\sU_1}
\mathop{\sum\sum}_{\substack{\xi_1,\xi_2\bmod{4^k}\\2\nmid \xi_1\xi_2}}\mathop{\sum\sum}_{\substack{u_1\sim U_1,u_2\sim U_2\\u_1\equiv\xi_1,u_2\equiv\xi_2\bmod{4^k}\\(2,u_1u_2)=(u_1,u_2)=1}}
\frac{Y_3-Y_2}{4^ku_1^2u_2^2}\\
&\ \ \ \ +2(1-2\delta)\sum_{k=0}^{k_0}\sum_{\xi\in\cR(2^k)}\sum_{\bU\in\sU_2}
\mathop{\sum\sum}_{\substack{\xi_1,\xi_2\bmod{4^k}\\2\nmid \xi_1\xi_2}}\mathop{\sum\sum}_{\substack{u_1\sim U_1,u_2\sim U_2\\u_1\equiv\xi_1,u_2\equiv\xi_2\bmod{4^k}\\(2,u_1u_2)=(u_1,u_2)=1\\ u_1\text{ is $U_1^{\theta}$-smooth}}}
\frac{Y_3-Y_2}{4^ku_1^2u_2^2}+O(x^{\frac{1}{2}}),
\end{align*}
where  $U_1,U_2$ are also restricted to be powers of 2.
Recall that
\[Y_3-Y_2=Y_3(2^ku_1u_2,\alpha)-Y_2(2^ku_1u_2)=2^ku_1u_2\sqrt{x}+O((2^ku_1u_2)^{1+\frac{1}{2\alpha}}).\]
We then obtain the lower bound
\begin{align*}
B(x,\alpha)&\geqslant2(1-2\delta)\Big(\sum_{k=0}^{k_0}\frac{\gamma(2^k)}{2^k}\Big)
\Bigg\{\mathop{\sum\sum}_{\substack{(u_1,u_2)\in\sV_1\\(2,u_1u_2)=(u_1,u_2)=1}}
\frac{1}{u_1u_2}+\mathop{\sum\sum}_{\substack{(u_1,u_2)\in\sV_2\\(2,u_1u_2)=(u_1,u_2)=1\\ u_1\text{ is $U_1^{\theta}$-smooth}}}
\frac{1}{u_1u_2}\Bigg\}\\
&\ \ \ \ \ +O(x^{\frac{1}{2}}\log x),
\end{align*}
where
\begin{align*}
\sV_1:=\{(u_1,u_2):u_1<u_2,\ \sqrt{x}< u_1u_2\leqslant x^\alpha,\  u_1\leqslant x^{\frac{1}{4}-\eta},\ u_2\leqslant u_1 x^{\frac{1}{2}-\eta}\},
\end{align*}
and
\begin{align*}
\sV_2:=\{(u_1,u_2):u_1<u_2,\ \sqrt{x}< u_1u_2\leqslant x^\alpha,\  u_1^{\frac{8+2\theta}{3}}u_2\leqslant x^{1-\eta},\ u_2\leqslant u_1 x^{\frac{1}{2}-\eta}\}\setminus \sV_1,
\end{align*}

Taking $k_0=k_0(\delta)$ very large, $\eta=\eta(\delta)$ very small, and letting $\delta$ tend to zero, we conclude from Lemma \ref{lm:basicasymptotics} that
\begin{align}\label{eq:B(x,alpha)-lowerbound-weak}
B(x,\alpha)&\geqslant\frac{1}{\pi^2}((2\alpha-1)(3-2\alpha)-o(1))\sqrt{x}\log^2x+B'(x,\alpha)
\end{align}
uniformly for $\alpha\in[\frac{1}{2},1],$  where
\begin{align*}
B'(x,\alpha)&=4\sqrt{x}\sum_{\substack{x^{\frac{1}{4}}<u_1\leqslant x^{\min\{\frac{\alpha}{2},\frac{3}{11+2\theta}\}}\\ u_1\text{ is $U_1^{\theta}$-smooth}\\2\nmid u_1}}\frac{\varphi(u_1)}{u_1^2}\log(\sqrt{x}u_1^{-\frac{5+2\theta}{3}}).
\end{align*}
Note that $B'(x,\alpha)$ is what we have gained more than Fouvry \cite{Fo16}. 
From Lemma \ref{lm:basicasymptotics-smoothnumbers}, we arrive at
\begin{align*}
B'(x,\alpha)&=\frac{1}{\pi^2}(F(\alpha,\theta)+o(1))\sqrt{x}\log^2x
\end{align*}
with 
\[
F(\alpha,\theta)=\rho\Big(\frac{1}{\theta}\Big)\Big\{4\min\Big\{\alpha,\frac{6}{11+2\theta}\Big\}-\frac{(10+4\theta)}{3}\min\Big\{\alpha,\frac{6}{11+2\theta}\Big\}^2-\frac{(7-2\theta)}{6}\Big\}.
\]
One may check that
\[F(\alpha,\theta)=\frac{1}{6}\rho\Big(\frac{1}{\theta}\Big)F_\theta(\alpha)\]
as given in Theorem \ref{thm:main}.
Combining this asymptotic evaluation for $B'(x,\alpha)$ with \eqref{eq:B(x,alpha)-lowerbound-weak}, we may conclude
a lower bound for $B(x,\alpha)$, from which and \eqref{eq:A(x,alpha)-asymptotic}, \eqref{eq:S=A+B}, we get
\begin{align}\label{eq:S(x,alpha)-lowerbound-weak}
S(x,\alpha)&\geqslant\frac{1}{\pi^2}\Big\{1+(2\alpha-1)(3-2\alpha)+\frac{1}{6}\rho\Big(\frac{1}{\theta}\Big)F_\theta(\alpha)-o(1)\Big\}\sqrt{x}\log^2x
\end{align}
uniformly for $\alpha\in[\frac{1}{2},1].$

\subsection{Concluding Theorem \ref{thm:main}}
To pass from a lower bound of $S(x,\alpha)$ to that of $S^\rf(x,\alpha)$, 
it is natural to invoke the identity \eqref{eq:S(x,alpha)-Sf(x,alpha)} and  Theorem \ref{thm:Hooley}. In fact, one can do a bit better following the arguments of Fouvry \cite{Fo16} and show that the above lower bound \eqref{eq:S(x,alpha)-lowerbound-weak}
also hold for $S^\rf(x,\alpha).$ This will depend on a more elaborate study of the contribution
from non-fundamental solutions. In other words, we would like to show that the non-fundamental
solutions create negligible contributions to $A(x,\alpha)$  and $B(x,\alpha).$

The following lemma is borrowed directly from \cite[Lemma 9.1]{Fo16}.
\begin{lemma}\label{lm:A-Af}
Uniformly for $\alpha\in[\frac{1}{2},1]$ and $x\geqslant2$, one has
\[A(x,\alpha)-A^\rf(x,\alpha)\ll \sqrt{x}\log x.\]
\end{lemma}

To deal with the contribution of the non-fundamental solutions to $B(x,\alpha)$, we also follow Fouvry. The above arguments which lead to \eqref{eq:S(x,alpha)-lowerbound-weak}
are essentially counting the number $\sN(x,\alpha;\varepsilon,k_0)$ of of 5-tuples of integers $(k,t,u_1,u_2,D)$ satisfying
\begin{align*}
2\nmid u_1u_2,\ \ (u_1,u_2)=1,\ \ 0\leqslant k\leqslant k_0,\ \ D\leqslant x,\end{align*}
\[X_{\frac{1}{2}}\leqslant 2^ku_1u_2\leqslant X_\alpha,\ \ t+2^ku_1u_2\sqrt{D}\leqslant D^{\frac{1}{2}+\alpha},\]
\[t^2\equiv1\bmod{4^k},\ \ t\equiv1\bmod{u_1^2},\ \ t\equiv-1\bmod{u_2^2},\]
as well as one of the following restrictions:
\begin{itemize}
\item $u_1\leqslant u_2\leqslant u_1x^{\frac{1}{2}-\eta},\ \ u_1\leqslant x^{\frac{1}{4}-\eta}$;
\item $u_2\leqslant u_1\leqslant u_2x^{\frac{1}{2}-\eta},\ \ u_2\leqslant x^{\frac{1}{4}-\eta}$;
\item $u_1\leqslant u_2\leqslant u_1x^{\frac{1}{2}-\eta},\ \ u_1^{\frac{8+2\theta}{3}}u_2\leqslant x^{1-\eta}$;
\item $u_2\leqslant u_1\leqslant u_2x^{\frac{1}{2}-\eta},\ \ u_1u_2^{\frac{8+2\theta}{3}}\leqslant x^{1-\eta}$.
\end{itemize}
By introducing the extra constraint $t+2^ku_1u_2\sqrt{D}=\varepsilon_D$,
we may also define $\sN^\rf(x,\alpha;\varepsilon,k_0)$. 
In fact, the above arguments yield
\begin{align}\label{eq:BN-BfNf}
B(x,\alpha)\geqslant \sN(x,\alpha;\varepsilon,k_0),\ \ B^\rf(x,\alpha)\geqslant \sN^\rf(x,\alpha;\varepsilon,k_0),
\end{align}
which are true for every positive $\varepsilon>0$ and for every $k_0\geqslant0.$
More precisely, we have proved for every $\delta$, $0<\varepsilon<\varepsilon_0(\delta),k_0>k_0(\delta)$ and $x>x_0(\delta)$ that
\begin{align}\label{eq:N(x,alpha)-lowerbound}
\sN(x,\alpha;\varepsilon,k_0)&\geqslant\frac{1}{\pi^2}\Big\{(2\alpha-1)(3-2\alpha)+\frac{1}{6}\rho\Big(\frac{1}{\theta}\Big)F_\theta(\alpha)-\delta\Big\}\sqrt{x}\log^2x
\end{align}
with $\alpha\in[\frac{1}{2},1].$

Following the approach of Fouvry \cite{Fo16}, we can state without proof that
\begin{lemma}\label{lm:N-Nf}
For every $k_0\geqslant0$ and every $\varepsilon>0$, one has
\[\sN(x,\alpha;\varepsilon,k_0)-\sN^\rf(x,\alpha;\varepsilon,k_0)\ll_{\varepsilon,k_0} \sqrt{x}\log x\]
uniformly for $\alpha\in[\frac{1}{2},1]$ and $x\geqslant2.$
\end{lemma}

We are now ready to complete the proof of Theorem \ref{thm:main}.
In view of \eqref{eq:BN-BfNf}, we may write
\begin{align*}
S^\rf(x,\alpha)
&=A^\rf(x,\alpha)+B^\rf(x,\alpha)\\
&\geqslant A^\rf(x,\alpha)+\sN^\rf(x,\alpha;\varepsilon,k_0)\\
&=A(x,\alpha)+\sN(x,\alpha;\varepsilon,k_0)+\{A^\rf(x,\alpha)-A(x,\alpha)\}\\
&\ \ \ \ \ +\{\sN^\rf(x,\alpha;\varepsilon,k_0)-\sN(x,\alpha;\varepsilon,k_0)\},
\end{align*}
which are true for every $\varepsilon>0$ and $k_0\geqslant0$. From Lemmas \ref{lm:A-Af} and \ref{lm:N-Nf}, we obtain
\begin{align*}
S^\rf(x,\alpha)
&\geqslant A(x,\alpha)+\sN(x,\alpha;\varepsilon,k_0)+O(\sqrt{x}\log x).
\end{align*}
By \eqref{eq:A(x,alpha)-asymptotic}, \eqref{eq:N(x,alpha)-lowerbound}, and by
choosing $k_0=k_0(\delta)$ sufficiently large, $\eta=\eta(\delta)$ sufficiently small, and letting $\delta$ tend to zero, we find the lower bound \eqref{eq:S(x,alpha)-lowerbound-weak}
holds definitely for $S^\rf(x,\alpha).$
This establishes \eqref{eq:Sf(x,alpha)-lowerbound}.

The lower bounds for $S(x,\alpha)$ in Theorem \ref{thm:main} can be deduced from \eqref{eq:Fouvry-Sf} by adding the contribution
of the non-fundamental solutions, as it is shown by \eqref{eq:S(x,alpha)-Sf(x,alpha)}.

\smallskip

\section{Estimate for triple exponential sums}\label{sec:expsums}

We now prove Theorem \ref{thm:tripleexpsum-fSestimate}. For the economy of the presentation, we only focus on the case $k=0$ and define
\begin{align}\label{eq:tripleexpsum-fS}
\fS(U_1,U_2,H):=\sum_{h\leqslant H}\mathop{\sum_{u_1\in]U_1,U_1^*]}\sum_{u_2\geqslant1}}
_{\substack{(2,u_1u_2)=(u_1,u_2)=1\\ u_1\text{ is $U_1^{\theta}$-smooth}}}
g_1\Big(\frac{u_2}{U_2}\Big)\ue\Big(\frac{h\overline{u_2^2}}{u_1^2}\Big).
\end{align}
We would like to show that
\begin{align}\label{eq:tripleexpsum-expected-smooth-simplifiedversion}
\fS(U_1,U_2,H)\ll U_1U_2x^{-\varepsilon_0},\ \ \ \ H=U_1U_2x^{-\frac{1}{2}+\varepsilon}
\end{align}
for some $\varepsilon_0\in~]0,10^{-2017}[$
while $U_1,U_2$ fall into the ranges in \eqref{eq:U1U2-sizes-smooth}.

By Poisson summation, the $u_2$-sum in \eqref{eq:tripleexpsum-fS} becomes
\begin{align*}
\sum_{u_2}&=\frac{U_2}{2u_1^2}\sum_{r\in\bZ}\widehat{g_1}\Big(\frac{rU_2}{2u_1^2}\Big)\sum_{\substack{z\bmod{2u_1^2}\\(z,2u_1)=1}}\ue\Big(\frac{h\overline{z^2}}{u_1^2}+\frac{rz}{2u_1^2}\Big).
\end{align*}
From the Chinese remainder theorem, the sum over $z$ can be rewritten as
\begin{align*}
&\ \ \ \ \sum_{\substack{z_1\bmod{u_1^2}\\(z_1,u_1)=1}}\sum_{\substack{z_2\bmod{2}\\(z,2)=1}}
\ue\Big(\frac{h\overline{(2z_1+u_1^2z_2)^2}}{u_1^2}+\frac{r(2z_1+u_1^2z_2)}{2u_1^2}\Big)\\
&=\sum_{\substack{z_1\bmod{u_1^2}\\(z_1,u_1)=1}}\sum_{\substack{z_2\bmod{2}\\(z,2)=1}}
\ue\Big(\frac{h\overline{(2z_1)^2}}{u_1^2}+\frac{rz_1}{u_1^2}+\frac{rz_2}{2}\Big)=(-1)^rK(r,\overline{4}h;u_1^2),
\end{align*}
where $K$ is an analogue of Kloosterman sums:
\begin{align*}
K(m,n;q)=\sideset{}{^*}\sum_{x\bmod q}\ue\Big(\frac{mx+n\overline{x}^2}{q}\Big).
\end{align*}
Hence we may conclude that
\begin{align*}
\sum_{u_2}&=\frac{U_2}{2u_1^2}\sum_{r\in\bZ}(-1)^r\widehat{g_1}\Big(\frac{rU_2}{2u_1^2}\Big)K(r,\overline{4}h;u_1^2).
\end{align*}

Note that $K(0,h;q)=T(h,q)$ and if $q$ is odd, we have $T(\overline{4}h,q)=T(h,q)$. 
According to $r=0$ and $r\neq0$, we split $\fS(U_1,U_2,H)$ by
\begin{align*}
\fS(U_1,U_2,H)&=\fS_1(U_1,U_2,H)+\fS_2(U_1,U_2,H),
\end{align*}
where
\begin{align*}
\fS_1(U_1,U_2,H)&=\frac{U_2(1-\delta)}{2}\sum_{h\leqslant H}\sum_{\substack{u_1\in]U_1,U_1^*]\\2\nmid u_1\text{ is $U_1^{\theta}$-smooth}}}
\frac{T(h,u_1^2)}{u_1^2}\end{align*}
and
\begin{align*}
\fS_2(U_1,U_2,H)&=\frac{U_2}{2}\sum_{h\leqslant H}\sum_{\substack{u_1\in]U_1,U_1^*]\\2\nmid u_1\text{ is $U_1^{\theta}$-smooth}}}
\frac{1}{u_1^2}\sum_{0\neq r\in\bZ}(-1)^r\widehat{g}_1\Big(\frac{rU_2}{2u_1^2}\Big)K(r,\overline{4}h;u_1^2).
\end{align*}
Following the approach of Fouvry, one may express $T(h,u_1^2)$
in terms of Jacobi symbols (see \cite[Lemma 6.2]{Fo16}) and then appeal to the bilinear estimate of Heath-Brown \cite{HB95}, getting
\begin{align}\label{eq:fS_1(U_1,U_2,H)-upperbound}
\fS_1(U_1,U_2,H)&\ll (U_1U_2H)^\varepsilon (HU_1U_2^{-1}+H^{\frac{1}{2}}U_1),
\end{align}
which produces the second restriction in \eqref{eq:U1U2-sizes-Fouvry}.
Moreover, Fouvry proved that
\begin{align*}
\fS_2(U_1,U_2,H)&\ll (U_1U_2H)^\varepsilon HU_2^2,
\end{align*}
which produces the first restriction in \eqref{eq:U1U2-sizes-Fouvry}.

Our task will be proving a stronger estimate for $\fS_2(U_1,U_2,H)$ by virtue of the special structure of $u_1$. More precisely, we shall prove that
\begin{align*}
\fS_2(U_1,U_2,H)&\ll (U_1U_2H)^\varepsilon HU_1U_2^{\frac{1}{2}}(Q+U_1^{\frac{1}{2}}Q^{-\frac{1}{2}}),
\end{align*}
where $Q$ will be chosen at our demand.
This would at least require the following inequality as proved by Fouvry \cite{Fo16}.
In fact, Fouvry only considered those $q$ that are perfect squares, and his argument also applies to more general $q$.
\begin{lemma}\label{lm:K(m,n;q)-upperbound}
Let $q$ be an odd positive integer. Then we have
\begin{align*}
|K(m,n;q)|\leqslant3^{\omega(q)} (m,n,q)\sqrt{q}.
\end{align*}
\end{lemma}

As an extension to $K(m,n;q)$, we define another exponential sum
\begin{align}\label{eq:B(m,n,ell,u;q)}
\sB(m,n,\ell,u;q)&=\sideset{}{^*}\sum_{a\bmod q}
\ue\Big(\frac{m\overline{a^2}+n\overline{(a+u)^2}+\ell a}{q}\Big),
\end{align}
where $*$ in summation reminds us to sum over primitive elements, i.e., $(a(a+u),q)=1.$
We will need the following inequality.
\begin{lemma}\label{lm:B(m,n,ell,u;q)-estimate}
For each odd $q\geqslant1$, we have
\begin{align*}
|\sB(m,n,\ell,u;q)|
&\leqslant 3q^{\frac{1}{2}}(\ell,m,n,q^\flat)^{\frac{1}{2}}(\ell,m+n,u,q^\flat)^{\frac{1}{2}}(\ell,m,n,q^\ddagger)(\ell,m+n,u,q^\ddagger)\cdot18^{\omega(q)}\cdot\Xi(q)^{\frac{1}{2}},
\end{align*}
where 
\[q^\ddagger=\prod_{p^2\| q}p,\ \ \ \ \Xi(q)= \prod_{p^\nu\|q, \, \nu\geqslant3} p^\nu.\]
\end{lemma}
The proof of Lemma \ref{lm:B(m,n,ell,u;q)-estimate} will be given in the last appendix.

We now start to prove Theorem \ref{thm:tripleexpsum-fSestimate}.
Due to the decay of $\widehat{g}_1$, we may truncate the $r$-sum in $\fS_2(U_1,U_2,H)$ by $R=U_1^{2+\varepsilon}U_2^{-1},$ so that
\begin{align}\label{eq:fS_1(U_1,U_2,H)-truncation}
\fS_2(U_1,U_2,H)&=\frac{U_2}{2}\sum_{h\leqslant H}\sum_{\substack{u_1\in]U_1,U_1^*]\\2\nmid u_1\text{ is $U_1^{\theta}$-smooth}}}
\frac{1}{u_1^2}\sum_{1\leqslant|r|\leqslant R}(-1)^r\widehat{g}_1\Big(\frac{rU_2}{2u_1^2}\Big)K(r,\overline{4}h;u_1^2)\\
&\ \ \ \ \ \ +O((U_1U_2H)^{-2017}).\nonumber
\end{align}
Our project will be controlling the cancellations while summing over $r$ with 
$1\leqslant|r|\leqslant R$. More precisely, we would like to estimate
\begin{align*}
\sum_{1\leqslant|r|\leqslant R}(-1)^r\widehat{g}_1\Big(\frac{rU_2}{2u_1^2}\Big)K(r,\overline{4}h;u_1^2).
\end{align*}
The contributions from negative $r$ and positive $r$ can be treated similarly, it thus suffices to consider
\begin{align*}
\sum_{r\leqslant R}(-1)^r\widehat{g}_1\Big(\frac{rU_2}{2u_1^2}\Big)K(r,\overline{4}h;u_1^2),
\end{align*}
which can be rewritten as
\begin{align*}
2\sum_{r\leqslant R/2}\widehat{g}_1\Big(\frac{rU_2}{u_1^2}\Big)K(2r,\overline{4}h;u_1^2)-\sum_{r\leqslant R}\widehat{g}_1\Big(\frac{rU_2}{2u_1^2}\Big)K(r,\overline{4}h;u_1^2).
\end{align*}
Clearly, $K(2r,\overline{4}h;u_1^2)=K(r,h;u_1^2)$.
By partial summation, it suffices to consider
\begin{align*}
\varSigma(R,u_1^2;h):=\sum_{r\leqslant R}K(r,h;u_1^2).
\end{align*}
For each $u_0\mid u_1$, we define $d,q_1,q_2$ by 
\[d=(u_0,(u_1/u_0)^\infty),\ \ q_1=(u_0/d)^2, \ \ q_2=u_1^2/q_1.\]
It follows that $u_1^2=q_1q_2$ and $(q_1,q_2)=1$. 
By virtue of the $q$-analogue of van der Corput method, we will prove
\begin{lemma}\label{lm:Sigma-upperbound}
With the above notation, we have
\begin{align*}
\varSigma(R,u_1^2;h)
&\ll u_1^{\varepsilon}R^{\frac{1}{2}}\{u_0u_1+u_0R^{\frac{1}{2}}(u_0,u_1/u_0)+(u_0,(u_1/u_0)^\infty)^{\frac{1}{2}}(h,u_0^2)^{\frac{1}{2}}u_1^{\frac{3}{2}}(u_1^\sharp)^{\frac{1}{4}}/u_0^{\frac{1}{2}}\}.\end{align*}
\end{lemma}

\proof Lemma \ref{lm:Sigma-upperbound} is a trivial consequence of Lemma \ref{lm:K(m,n;q)-upperbound} if $R\leqslant u_0^2.$ We now assume $R> u_0^2.$
Denote by $I_R$ the characteristic function of the interval $[1,R].$ Thus,
\begin{align*}
\varSigma(R,u_1^2;h)=\sum_{r\in\bZ}I_R(r)K(r,h;u_1^2)=\sum_{r\in\bZ}I_R(r+u_0^2\ell)
K(r+u_0^2\ell,h;u_1^2)\end{align*}
for any $\ell\in\bZ.$ From the Chinese remainder theorem, we may write
\begin{align*}
K(r+u_0^2\ell,h;u_1^2)
&=K(r+u_0^2\ell,h\overline{q}_2^3;q_1)K(r+u_0^2\ell,h\overline{q}_1^3;q_2)\\
&=K(r,h\overline{q}_2^3;q_1)K(r+u_0^2\ell,h\overline{q}_1^3;q_2),
\end{align*}
where we have used the fact that $q_1\mid u_0^2.$

For $L=[R/u_0^2]$, we sum over $\ell$, getting
\begin{align*}
\varSigma(R,u_1^2;h)&=\frac{1}{L}\sum_{\ell\leqslant L}\sum_{r\in\bZ}I_R(r+u_0^2\ell)
K(r+u_0^2\ell,h;u_1^2)\\
&\leqslant\frac{1}{L}\sum_{r\in\bZ}|K(r,h\overline{q}_2^3;q_1)|\Bigg|\sum_{\ell\leqslant L}I_R(r+u_0^2\ell)
K(r+u_0^2\ell,h\overline{q}_1^3;q_2)\Bigg|.\end{align*}
From Lemma \ref{lm:K(m,n;q)-upperbound} we find
\begin{align*}
\varSigma(R,u_1^2;h)
&\ll\frac{q^\varepsilon\sqrt{q_1}}{L}\sum_{r\in\bZ}(r,h,q_1)|\Bigg|\sum_{\ell\leqslant L}I_R(r+u_0^2\ell)
K(r+u_0^2\ell,h\overline{q}_1^3;q_2)\Bigg|.\end{align*}
In view of the support of $I_R$, the sum over $r$ is in fact restricted to $[-R,R].$
By Cauchy inequality, we derive that
\begin{align*}
\varSigma(R,u_1^2;h)^2
&\ll\frac{q^\varepsilon q_1R}{L^2}\sum_{r\in\bZ}\Bigg|\sum_{\ell\leqslant L}I_R(r+u_0^2\ell)
K(r+u_0^2\ell,h\overline{q}_1^3;q_2)\Bigg|^2.\end{align*}
Squaring out and switching summations, we get
\begin{align}\label{eq:Sigma^2-upperbound}
\varSigma(R,u_1^2;h)^2
&\ll\frac{q^\varepsilon q_1R}{L}\sum_{0\leqslant|\ell|\leqslant L}\Bigg|\sum_{r\in I_\ell}K(r,h\overline{q}_1^3;q_2)\overline{K(r+u_0^2\ell,h\overline{q}_1^3;q_2)}\Bigg|,\end{align}
where $I_\ell$ is an interval, depending on $\ell$, of length at most $R$.
For $\ell=0$, we appeal to Lemma \ref{lm:K(m,n;q)-upperbound} to estimate the $r$-sum trivially. For $1\leqslant|\ell|\leqslant L$,
reasonable cancellations in the $r$-sum are expected. In fact, by completion, we have
\begin{align*}
\sum_{r\in I_\ell}[\cdots]&=\sum_{y\bmod{q_2}}I(y,q_2)\widehat{K}(y,q_2;h,q_1,u_0,\ell),\end{align*}
where
\begin{align*}
I(y,q_2)&=\sum_{r\in I_\ell}\ue\Big(\frac{-yr}{q_2}\Big),\end{align*}
and
\begin{align*}
\widehat{K}(y,q_2;h,q_1,u_0,\ell)&=\frac{1}{q_2}\sum_{x\bmod{q_2}}K(x,h\overline{q}_1^3;q_2)\overline{K(x+u_0^2\ell,h\overline{q}_1^3;q_2)}
\ue\Big(\frac{yx}{q_2}\Big).\end{align*}
On one hand,
\begin{align*}
I(y,q_2)&\ll\min\Big\{R,\Big\|\frac{y}{q_2}\Big\|^{-1}\Big\}.\end{align*}
Opening each $K$ by definition, the orthogonality of additive characters gives
\begin{align*}
\widehat{K}(y,q_2;h,q_1,u_0,\ell)&=\ue\Big(\frac{-u_0^2\ell y}{q_2}\Big)\sideset{}{^*}\sum_{a\bmod{q_2}}
\ue\Big(\frac{h\overline{q}_1^3(\overline{a^2}-\overline{(a+y)^2})-u_0^2\ell a}{q_2}\Big)\\
&=\ue\Big(\frac{-u_0^2\ell y}{q_2}\Big)\sB(h\overline{q}_1^3,-h\overline{q}_1^3,-u_0^2\ell,y;q_2).\end{align*}
For $y=0$, we have
\[\sB(h\overline{q}_1^3,-h\overline{q}_1^3,-u_0^2\ell,0;q_2)=\sideset{}{^*}\sum_{a\bmod{q_2}}
\ue\Big(\frac{u_0^2\ell a}{q_2}\Big)\ll(u_0^2\ell,q_2).\]
For $y\neq0$, we would like to appeal to Lemma \ref{lm:B(m,n,ell,u;q)-estimate}.
To do so, we first derive from Lemma \ref{lm:B(m,n,ell,u;q)-estimate} that
(we have $q_2^\flat=1$ since $q_2$ is a perfect square)
\begin{align*}
|\sB(h\overline{q}_1^3,-h\overline{q}_1^3,-u_0^2\ell,y;q_2)|\ll
q_2^{\frac{1}{2}+\varepsilon}(h,u_0^2\ell,q_2^\ddagger)(y,u_0^2\ell,q_2^\ddagger)\cdot\Xi(q_2)^{\frac{1}{2}}.\end{align*}
This yields
\begin{align*}
\sum_{r\in I_\ell}[\cdots]&\ll R(u_0^2\ell,q_2)+q_2^{\frac{1}{2}+\varepsilon}(h,u_0^2\ell,q_2^\ddagger)\cdot\Xi(q_2)^{\frac{1}{2}}\sum_{1\leqslant|y|\leqslant q_2/2}\min\Big\{R,\frac{q_2}{y}\Big\}(y,u_0^2\ell,q_2^\ddagger)\\
&\ll R(u_0^2\ell,q_2)+q_2^{\frac{3}{2}+\varepsilon}(h,u_0^2\ell,q_2^\ddagger)\cdot\Xi(q_2)^{\frac{1}{2}},\end{align*}
from which and \eqref{eq:Sigma^2-upperbound} we obtain 
\begin{align*}
\varSigma(R,u_1^2;h)^2
&\ll\frac{u_1^\varepsilon q_1R}{L}\{Rq_2(h,q_2)+LR(u_0^2,q_2)+Lq_2^{\frac{3}{2}}(h,u_0^2,q_2^\ddagger)\cdot\Xi(q_2)^{\frac{1}{2}}\}\\
&\ll u_1^{\varepsilon}R\{(u_0u_1)^2+q_1R(u_0^2,q_2)+u_1^2q_2^{\frac{1}{2}}(h,u_0^2)\cdot\Xi(q_2)^{\frac{1}{2}}\}.\end{align*}
Note that $q_2=(u_1d/u_0)^2=\{u_1(u_0,(u_1/u_0)^\infty)/u_0\}^2.$ Thus,
\begin{align*}
\Xi(q_2)=\Xi(\{u_1(u_0,(u_1/u_0)^\infty)/u_0\}^2)=(u_1(u_0,(u_1/u_0)^\infty)/u_0)^\sharp
\leqslant u_1^\sharp.\end{align*}
We then conclude that
\begin{align*}
\varSigma(R,u_1^2;h)^2
&\ll u_1^{\varepsilon}R\{(u_0u_1)^2+u_0^2R(u_0,u_1/u_0)^2+(u_0,(u_1/u_0)^\infty)(h,u_0^2)u_1^3(u_1^\sharp)^{\frac{1}{2}}/u_0\},\end{align*}
which gives Lemma \ref{lm:Sigma-upperbound} immediately.
\endproof

In view of Lemma \ref{lm:Sigma-upperbound} and the discussions before it, we may derive from \eqref{eq:fS_1(U_1,U_2,H)-truncation} that
\begin{align*}
\fS_2(U_1,U_2,H)&\ll U_1^\varepsilon HQU_2(R^{\frac{1}{2}}+RU_1^{-1})+\frac{U_1^\varepsilon HU_2R^{\frac{1}{2}}}{(QU_1)^{\frac{1}{2}}}\sum_{u_1\sim U_1}\sum_{ab=u_1}
(a,b^\infty)^{\frac{1}{2}}((ab)^{\sharp})^{\frac{1}{4}}.
\end{align*}
where $U_1^{\theta_0}\leqslant Q\leqslant U_1^{\theta_0+\theta}$. In view of the choice $R=U_1^{2+\varepsilon}U_2^{-1}$, Lemma \ref{lm:Rankin} yields
\begin{align*}
\fS_2(U_1,U_2,H)
&\ll U_1^\varepsilon HU_1U_2^{\frac{1}{2}}\{Q+U_1^{\frac{1}{2}}Q^{-\frac{1}{2}}\},
\end{align*}
from which and \eqref{eq:fS_1(U_1,U_2,H)-upperbound} we conclude that
\begin{align*}
\fS(U_1,U_2,H)&\ll (U_1U_2H)^\varepsilon \{HU_1U_2^{-1}+H^{\frac{1}{2}}U_1+HU_1U_2^{\frac{1}{2}}(Q+U_1^{\frac{1}{2}}Q^{-\frac{1}{2}})\}.
\end{align*}
Choosing $\theta_0=\frac{1-2\theta}{3}$, we then have
 $U_1^{\frac{1-2\theta}{3}}\leqslant Q\leqslant U_1^{\frac{1+\theta}{3}}$. Recalling the choice of $H$, we then arrive at the expected estimate \eqref{eq:tripleexpsum-expected-smooth-simplifiedversion}, provided that
\eqref{eq:U1U2-sizes-smooth} holds.

\smallskip

\appendix
\section{Mean values of arithmetic functions}
\subsection{Some basic asymptotics}
The first part of the appendix is devoted to state several basic asymptotics.

Recall that $\gamma(u)$ denotes the number of solutions to the congruence equation $x^2\equiv1\bmod{u^2}.$ As a multiplicative function, $\gamma$ satisfies the evaluation \eqref{eq:gamma-primepower}. We are now ready to state the following averages.

\begin{lemma}\label{lm:basicasymptotics}
As $N\rightarrow+\infty$, we have
\begin{align*}
\sum_{0\leqslant n\leqslant N}\frac{\gamma(2^n)}{2^n}&=4+O(2^{-N}),
\end{align*}
\begin{align*}
\sum_{\substack{n\leqslant N\\(n,q)=1}}\frac{1}{n}=\frac{\varphi(q)}{q}\log N+O\Big(1+\sum_{d|q}\frac{\log d}{d}\Big),
\end{align*}
\begin{align*}
\sum_{\substack{n\leqslant N\\(n,2)=1}}\frac{\varphi(n)}{n^2}\sim\frac{4}{\pi^2}\log N,
\end{align*}
and
\begin{align*}
\sum_{\substack{n\leqslant N\\(n,2)=1}}\frac{\varphi(n)\log n}{n^2}\sim\frac{2}{\pi^2}\log^2N.
\end{align*}
\end{lemma}

\proof
The first one can be derived from the evaluations of $\gamma(2^n)$ as given in \eqref{eq:gamma-primepower}. The other three asymptotics can be found
in Fouvry \cite[Lemma 8.1, Lemma 8.2]{Fo16}.
\endproof

\subsection{Smooth numbers}
Denote by $\cS(x,y)$ the set of $y$-smooth numbers not exceeding $x$. Write
$\Psi(x,y)=|\cS(x,y)|.$
We now introduce the Dickman function $\rho(u)$ by
\begin{align}\label{eq:Dickman}
\begin{cases}
\rho(u)=1,\ \ \ &u\in~]0,1],\\
u\rho'(u)+\rho(u-1)=0,&u\in~]1,+\infty[.
\end{cases}
\end{align}
In the first several intervals, we have
\begin{align*}
\rho(u)&=1-\log u,\ \ \ u\in~]1,2],\\
\rho(u)&=1-\log u+\int_2^u\frac{\log(t-1)}{t}\ud t,\ \ \ u\in~]2,3],\\
\rho(u)&=1-\log u+\int_2^u\frac{\log(t-1)}{t}\ud t-\int_3^u\frac{\ud t}{t}\int_2^{t-1}\frac{\log(s-1)}{s}\ud s,\ \ \ u\in~]3,4].
\end{align*}

The following lemma is classical and shows $\rho$ is the density function of smooth numbers.
\begin{lemma}\label{lm:smoothnumbers}
Uniformly for $x\geqslant y\geqslant2$, we have
\begin{align*}
\Psi(x,y)=x\rho\Big(\frac{\log x}{\log y}\Big)+O\Big(\frac{x}{\log y}\Big).
\end{align*}
\end{lemma}

\proof See \cite[P. 367, Theorem 6]{Te95}.
\endproof

\begin{lemma}\label{lm:basicasymptotics-smoothnumbers}
Let $\theta\in~]0,1[$ be fixed. As $N\rightarrow+\infty$, we have
\begin{align*}
\sum_{\substack{2\nmid n\leqslant N\\n\text{\rm~ is $n^{\theta}$-smooth}}}\frac{\varphi(n)}{n^2}\sim\frac{4}{\pi^2}\rho\Big(\frac{1}{\theta}\Big)\log N,
\end{align*}
and
\begin{align*}
\sum_{\substack{2\nmid n\leqslant N\\n\text{\rm~ is $n^{\theta}$-smooth}}}\frac{\varphi(n)\log n}{n^2}\sim\frac{2}{\pi^2}\rho\Big(\frac{1}{\theta}\Big)\log^2N.
\end{align*}
\end{lemma}

\proof
One can refer to \cite{TW03}, for instance, for some general theorems on the mean values of multiplicative functions over smooth numbers. In particular, one has
\begin{align*}
\sum_{\substack{2\nmid n\leqslant N\\n\text{\rm~ is $n^{\theta}$-smooth}}}\frac{\varphi(n)}{n}=\frac{4}{\pi^2}\rho\Big(\frac{1}{\theta}\Big) N+O\Big(\frac{N}{\log N}\Big).
\end{align*}
The lemma then follows from the partial summation.
\endproof

\begin{lemma}\label{lm:squarefullpart-average}
For all $q\geqslant1$, we have
\begin{align*}
\sum_{n\leqslant x}(n^{\sharp})^{\frac{1}{4}}(n,q)^{\frac{1}{2}}\ll x(qx)^\varepsilon.
\end{align*}
\end{lemma}
\proof
Note that $n\mapsto (n^{\sharp})^{\frac{1}{4}}(n,q)^{\frac{1}{2}}$ is multiplicative. For $\Re s>1,$ we consider the Dirichet series
\begin{align*}
\sD(s)=\sum_{n\geqslant1}\frac{(n^{\sharp})^{\frac{1}{4}}(n,q)^{\frac{1}{2}}}{n^s}.
\end{align*}
By Euler product formula, we have
\begin{align*}
\sD(s)&=\prod_p\Big(\sum_{k\geqslant0}\frac{((p^k)^{\sharp})^{\frac{1}{4}}(p^k,q)^{\frac{1}{2}}}{p^{ks}}\Big)\\
&=\prod_p\Big(1+\frac{(p,q)^{\frac{1}{2}}}{p^{s}}+\sum_{k\geqslant2}\frac{(p^k,q)^{\frac{1}{2}}}{p^{k(s-1/4)}}\Big)\\
&=\prod_{p\nmid q}\Big(1+\frac{1}{p^{s}}+\sum_{k\geqslant2}\frac{1}{p^{k(s-1/4)}}\Big)\cdot \prod_{p\mid q}\Big(1+\frac{1}{p^{s-1/2}}+\sum_{k\geqslant2}\frac{(p^k,q)^{\frac{1}{2}}}{p^{k(s-1/4)}}\Big)\\
&=\zeta(s)\sD_1(s)\sD_2(s,q),
\end{align*}
where $\sD_1(s)$ is holomorphic for $\Re s>3/4$ and
\begin{align*}
\sD_2(s,q)&=\prod_{p\mid q}\Big(1+\frac{1}{p^{s}}+\sum_{k\geqslant2}\frac{1}{p^{k(s-1/4)}}\Big)^{-1}\cdot \Big(1+\frac{1}{p^{s-1/2}}+\sum_{k\geqslant2}\frac{(p^k,q)^{\frac{1}{2}}}{p^{k(s-1/4)}}\Big).
\end{align*}
The lemma then follows from a routine application of Perron's formula.
\endproof

The following inequality is a consequence of Rankin's method, which is a stronger version of \cite[Lemma 7.2]{Fo16}.
\begin{lemma}\label{lm:Rankin}
For any $\varepsilon>0$, one has
\begin{align*}
\sum_{mn\leqslant N}(m,n^\infty)^{\frac{1}{2}}((mn)^{\sharp})^{\frac{1}{4}}&\ll N^{1+\varepsilon}.
\end{align*}
\end{lemma}
\proof
Denote by $S$ the sum in question.
First, we have
\begin{align*}
S&\ll\sum_{n\leqslant N}\sum_{d\mid n^\infty}d^{\frac{1}{2}}\sum_{\substack{m\leqslant N/n\\d\mid m}}((mn)^{\sharp})^{\frac{1}{4}}
=\sum_{n\leqslant N}\sum_{\substack{d\leqslant N\\ d\mid n^\infty}}d^{\frac{1}{2}}\sum_{\substack{m\leqslant N/(nd)}}((mnd)^{\sharp})^{\frac{1}{4}}.
\end{align*}
Note that
\begin{align*}
(ab)^{\flat}&=\prod_{p\| ab}p=\prod_{p\| a,p\nmid b}p\cdot \prod_{p\| b,p\nmid a}p
=a^\flat\cdot b^\flat\cdot \prod_{p\| a,p\mid b}\frac{1}{p}\cdot \prod_{p\| b,p\mid a}\frac{1}{p}\geqslant\frac{a^\flat\cdot b^\flat}{\{(a,b)^\flat\}^2},
\end{align*}
giving
\begin{align*}
(ab)^{\sharp}&\leqslant a^\sharp\cdot b^\sharp\cdot\{(a,b)^\flat\}^2,
\end{align*}
from which and Lemma \ref{lm:squarefullpart-average} it follows that
\begin{align*}
S&\ll\sum_{n\leqslant N}\sum_{\substack{d\leqslant N\\ d\mid n^\infty}}d^{\frac{1}{2}}((nd)^{\sharp})^{\frac{1}{4}}\sum_{\substack{m\leqslant N/(nd)}}(m^{\sharp})^{\frac{1}{4}}(m,nd)^{\frac{1}{2}}\\
&\ll N^{1+\varepsilon}\sum_{n\leqslant N}\frac{1}{n}\sum_{\substack{d\leqslant N\\ d\mid n^\infty}}d^{-\frac{1}{2}}((nd)^{\sharp})^{\frac{1}{4}}\\
&\ll N^{1+\varepsilon}\sum_{n\leqslant N}\frac{(n^{\sharp})^{\frac{1}{4}}}{n}\sum_{\substack{d\leqslant N\\ d\mid n^\infty}}d^{-\frac{1}{2}}(d^{\sharp})^{\frac{1}{4}}((n,d)^\flat)^{\frac{1}{2}}.
\end{align*}
By Rankin's method, the last sum over $d$ is, for any $\varepsilon>0$,
\begin{align*}
&\leqslant N^\varepsilon\sum_{\substack{d\leqslant N\\ d\mid n^\infty}}d^{-\frac{1}{2}-\varepsilon}(d^{\sharp})^{\frac{1}{4}}((n,d)^\flat)^{\frac{1}{2}}\\
&\leqslant N^\varepsilon\prod_{p\mid n}\Big(1+\sum_{k\geqslant1}\frac{((p^k)^{\sharp})^{\frac{1}{4}}((n,p^k)^\flat)^{\frac{1}{2}}}{p^{k(1/2+\varepsilon)}}\Big)\\
&=N^\varepsilon\prod_{p\mid n}\Big(1+\frac{2}{p^{\varepsilon}}+\sum_{k\geqslant3}\frac{1}{p^{k(1/4+\varepsilon)-1/2}}\Big)\\
&\ll (nN)^\varepsilon
\end{align*}
by re-defining $\varepsilon$. 
We now get
\begin{align*}
S
&\ll N^{1+\varepsilon}\sum_{n\leqslant N}\frac{(n^{\sharp})^{\frac{1}{4}}}{n}\ll N^{1+\varepsilon},
\end{align*}
where the last step follows from Lemma \ref{lm:squarefullpart-average} together with partial summation and taking $q=1$ therein.
\endproof

\section{Estimate for $\sB(m,n,\ell,u;q)$}

While invoking the ideas of $q$-analogue of van der Corput method, we have transformed the original algebraic exponential sum $K(m,n,q)$ to a new sum $\sB(m,n,\ell,u;q)$
as given by \eqref{eq:B(m,n,ell,u;q)}. This appendix will be devoted to present an estimation for $\sB(m,n,\ell,u;q)$ that suits well in our applications to Theorem \ref{thm:main}. In fact, the job can be done for the general 
complete exponential sum
\begin{align*}
\Sigma(\lambda,q)
:= \sum_{a\bmod q} \mathrm{e}\Big(\frac{\lambda(a)}{q}\Big).
\end{align*}
Here $\lambda=\lambda_1/\lambda_2$ with $\lambda_1,\lambda_2\in\bZ[X]$ 
and $\lambda_1,\lambda_2$ being coprime in $\bZ[X]$.
The values of $a$ such that $(\lambda_2(a),q)\neq1$ are excluded from summation. 
We define the degree of $\lambda$ by
\[d=d(\lambda)=\deg(\lambda_1)+\deg(\lambda_2).\]
If 
\[\lambda_1(x)=\sum_{0\leqslant j\leqslant d_1}r_jx^j\in\bZ[x],\ \ \ \lambda_2(x)=\sum_{0\leqslant j\leqslant d_2}t_jx^j\in\bZ[x],\]
we then adopt the convention that 
\begin{align}
(\lambda,c)_*&=(r_0,r_1,r_2,\cdots,r_{d_1},c),
\label{eq:gcdrationalfunctions1}\\
(\lambda,c)&=(r_1,r_2,\cdots,r_{d_1},t_1,t_2,\cdots,t_{d_2},c),\label{eq:gcdrationalfunctions2}\\ (\lambda',c)&=(\lambda_1'\lambda_2-\lambda_1\lambda_2',c)
\label{eq:gcdrationalfunctions3}
\end{align}
for all $c\mid q.$

There are many known estimates for complete exponential sums in the literature. In \cite[Theorem B.1]{WX16}, we obtained the following estimate for $\Sigma(\lambda,q)$.

\begin{theorem}\label{thm:S(lambda,q)bound}
Let $d=d(\lambda)\geqslant1.$ For $q\geqslant1,$ we have
\begin{align*}
|\Sigma(\lambda,q)|
\leqslant q^{\frac{1}{2}} (\lambda,q^\flat)^{\frac{1}{2}}(\lambda,q^\flat)_*^{\frac{1}{2}}(\lambda',q^\ddagger)(2d)^{\omega(q)}\cdot\Xi(q)^{\frac{1}{2}},
\end{align*}
where $q^\ddagger$ and $\Xi(q)$ are given as in Lemma {\rm \ref{lm:B(m,n,ell,u;q)-estimate}}.
\end{theorem}
In fact, the case of $q$ being a prime is essentially due to Weil \cite{We48}, and 
one can refer to \cite{Bom66} for a complete proof that suits quite well in our situation.
Thanks to the Chinese remainder theorem, the evaluation of 
$\Sigma(\lambda,q)$ can be reduced to the case of prime power moduli $p^\alpha (\alpha\geqslant2)$, which can usually be treated following an elementary device; see \cite[Section 12.3]{IK04} for details. This is in fact the line of the proof in \cite{WX16}.

The upper bound in Theorem \ref{thm:S(lambda,q)bound} is complicated at first glance and it even gives a worse bound than the trivial one if $\lambda$ is a constant function mod $q$.
However, 
it provides an essentially optimal bound if we have an extra average over $q$, since 
$(\lambda,q^\flat)^{\frac{1}{2}}(\lambda,q^\flat)_*^{\frac{1}{2}}(\lambda',q^\ddagger)\Xi(q)^{\frac{1}{2}}$ does not oscillate too much on average.

Given $m,n,u,\ell\in\bZ,$ we take
\begin{align*}
\lambda(x)=\frac{m}{x^2}+\frac{n}{(x+u)^2}+\ell x=\frac{\lambda_1(x)}{\lambda_2(x)}\end{align*}
with
\begin{align*}
\lambda_1(x)=\ell x^5+2\ell ux^4+\ell u^2x^3+(m+n)x^2+2mux+mu^2,\ \ \ \lambda_2(x)=x^2(x+u)^2.\end{align*}
Moreover, we have
\begin{align*}
&\ \ \ \ (\lambda_1'\lambda_2-\lambda_1\lambda_2')(x)\\
&=(5\ell x^4+8\ell ux^3+3\ell u^2x^2+2(m+n)x+2mu)x^2(x+u)^2\\
&\ \ \ \ -(4x^3+6ux^2+2u^2x)(\ell x^5+2\ell ux^4+\ell u^2x^3+(m+n)x^2+2mux+mu^2)\\
&=5\ell x^8+8\ell ux^7+3\ell u^2x^6+2(m+n)x^5+2mux^4\\
&\ \ \ \ +10\ell ux^7+16\ell u^2x^6+6\ell u^3x^5+4(m+n)ux^4+4mu^2x^3\\
&\ \ \ \ + 5\ell u^2x^6+8\ell u^3x^5+3\ell u^4x^4+2(m+n)u^2x^3+2mu^3x^2\\
&\ \ \ \ -(4\ell x^8+8\ell ux^7+4\ell u^2x^6+4(m+n)x^5+8mux^4+4mu^2x^3)\\
&\ \ \ \ -(6\ell ux^7+12\ell u^2x^6+6\ell u^3x^5+6(m+n)ux^4+12mu^2x^3+6mu^3x^2)\\
&\ \ \ \ -(2\ell u^2x^6+4\ell u^3x^5+2\ell u^4x^4+4mu^3x^2+2mu^4x)\\
&=\ell x^8+4\ell ux^7-6\ell u^2x^6-(2m+2n-4\ell u^3)x^5+(\ell u^4-2nu-8mu)x^4\\
&\ \ \ \ \ +(2n-10m)u^2x^3-8mu^3x^2-2mu^4x.\end{align*}
Thus, for odd $c\mid q$, we have
\begin{align*}
(\lambda,c)_*&=(\ell,2\ell u,\ell u^2,m+n,2mu,mu^2,c)=(\ell,m+n,mu,c)\leqslant(\ell,m,n,c)(\ell,m+n,u,c),\end{align*}
\begin{align*}
(\lambda,c)&=(\ell,2\ell u,\ell u^2,m+n,2mu,1,2u,u^2,c)=1\end{align*}
and
\begin{align*}
(\lambda',c)&=(\ell,4\ell u,6\ell u^2,2m+2n-4\ell u^3,\ell u^4-2nu-8mu,(2n-10m)u^2,8mu^3,2mu^4,c)\\
&=(\ell,m+n,(4m+n)u,(5m-n)u^2,mu^3,c)\\
&=(\ell,m+n,3mu,6mu^2,mu^3,c)\\
&\leqslant3(\ell,m+n,mu,c)\\
&\leqslant3(\ell,m,n,c)(\ell,m+n,u,c).\end{align*}
Note that $d(\lambda)=\deg(\lambda_1)+\deg(\lambda_2)=5+4=9.$

From Theorem \ref{thm:S(lambda,q)bound}, we may conclude the following estimate as given in Lemma \ref{lm:B(m,n,ell,u;q)-estimate}.
\begin{theorem}
For each odd $q\geqslant1$, we have
\begin{align*}
|\sB(m,n,\ell,u;q)|
&\leqslant 3q^{\frac{1}{2}}(\ell,m,n,q^\flat)^{\frac{1}{2}}(\ell,m+n,u,q^\flat)^{\frac{1}{2}}(\ell,m,n,q^\ddagger)(\ell,m+n,u,q^\ddagger)\cdot18^{\omega(q)}\cdot\Xi(q)^{\frac{1}{2}}.
\end{align*}
\end{theorem}

\endproof

\bigskip

\bibliographystyle{plainnat}

\end{document}